\documentclass[10pt,twoside]{article}
\usepackage[latin1]{inputenc}
\usepackage{amsmath}
\usepackage{graphicx}
\usepackage{yhmath}
\usepackage[table]{xcolor}
\usepackage{mathrsfs} 
\usepackage{amssymb}
\usepackage{url}
\usepackage{makecell}
\usepackage{arydshln}

\usepackage{stmaryrd} 

\input epsf
\def\limiten{\renewcommand{\arraystretch}{0.5}
\begin{array}[t]{c}\stackrel{}{\longrightarrow} \\
{\scriptstyle n\rightarrow
\infty}\end{array}\renewcommand{\arraystretch}{1}}

\def\limitepsn{\renewcommand{\arraystretch}{0.5}
\begin{array}[t]{c}\stackrel{a.s.}{\longrightarrow} \\
{\scriptstyle n \rightarrow
\infty}\end{array}\renewcommand{\arraystretch}{1}}

\def\limiteloin{\renewcommand{\arraystretch}{0.5}
\begin{array}[t]{c}\stackrel{{\cal D}}{\longrightarrow} \\
{\scriptstyle n\rightarrow
\infty}\end{array}\renewcommand{\arraystretch}{1}}

\def\limiteproban{\renewcommand{\arraystretch}{0.5}
\begin{array}[t]{c}\stackrel{{\cal P}}{\longrightarrow} \\
{\scriptstyle n\rightarrow
\infty}\end{array}\renewcommand{\arraystretch}{1}}

\numberwithin{equation}{section}

\newtheorem{thm}{Theorem}[section]

\newtheorem{example}[thm]{Exemple}

\newtheorem{lem}[thm]{Lemma}

\newtheorem{prop}[thm]{Proposition}

\newcommand{\E}{\ensuremath{\mathbb{E}}}
\newcommand{\R}{\ensuremath{\mathbb{R}}}
\newcommand{\Z}{\ensuremath{\mathbb{Z}}}

\newcommand{\N}{\ensuremath{\mathbb{N}}}

\newcommand{\var}{\ensuremath{\mathrm{Var}}}

\definecolor{grisclair}{gray}{0.9}

\renewcommand{\arraystretch}{.8}

\begin{document}
\title{\bf Consistent model selection procedure for general integer-valued time series }
 \maketitle \vspace{-1.0cm}
 \begin{center}
   Mamadou Lamine DIOP \footnote{Supported by the CY Advanced Studies (CY Cergy Paris Université, France), and the MME-DII center of excellence (ANR-11-LABEX-0023-01)} 
   and 
    William KENGNE  \footnote{Developed within the ANR BREAKRISK : ANR-17-CE26-0001-01} 
 \end{center}

  \begin{center}
 {\it  THEMA, CY Cergy Paris Université, 33 Boulevard du Port, 95011 Cergy-Pontoise Cedex, France.\\
  E-mail: mamadou-lamine.diop@u-cergy.fr ; william.kengne@u-cergy.fr  \\}
\end{center}

 \pagestyle{myheadings}
 \markboth{Consistent model selection procedure for general integer-valued time series}{Diop and Kengne}

\textbf{Abstract} :
 This paper deals with the problem of model selection for a general class of integer-valued time series. 
   We propose a penalized criterion based on the Poisson quasi-likelihood of the model. 
 Under certain regularity conditions, the consistency of the procedure as well as the consistency and the asymptotic normality of the Poisson quasi-likelihood estimator of the selected model are established.
 Simulation experiments are conducted for some classical models such as Poisson, binary INGARCH and  negative binomial model with nonlinear dynamic. Also, an application to a real dataset is provided.  


\section{Introduction}

Model selection is a fundamental step in many statistical analysis. There consists in choosing a model amongst
a collection such that its realizations are close as possible to that of reality manifested in the observed data. 
 Many authors have already examined this question by defining a reasonable setting related to the nature of the experiment or the data, which allows them to reduce the possible set of model that can be used.  
  Two approaches are often used for solving this question: the cross validation and the penalized contrast. 
  In cross validation, observations are split: part of data (the training sample) is used for fitting each competing model and  the remaining part (the validation sample) is used to measure the performances of models; and the model with the best overall performance is selected.
 %
It is a selection method based on the predictive ability of the models.
 See, for example, Stone (1974) and Allen (1974) for a pioneer works introduced.  
  In the procedures using a penalized contrast, the selection is done by minimizing a trade-off between a sum of an empirical risk which indicates how the model fits the data (estimation error) and a measure of model's complexity (approximation error). 
  The first examples of this kind of criterion are the Mallows' $C_p$ within the regression framework (see \cite{Mallows1974}) and the AIC (see \cite{Akaike1973}) or the BIC (see \cite{Schwarz1978}) criteria within the maximum-likelihood framework.   
We refer to Leeb and P\"{o}tscher (2009) for a relevant literature on these approaches. 

\medskip
     \medskip 
The question of interest here is about modelling time series of counts. 
These models have been a growing interest, given the large number of articles written in this direction during the last two decades; we refer to Fokianos {\it et al.} (2009), Doukhan {\it et al.} (2012), Fokianos and Neumann (2013), Doukhan and Kengne (2015), Davis and Liu (2016), Ahmad and Francq (2016) among others for some recent works.   
 The problem of model selection in such models has been addressed by several authors. For instance, Peng {\it et al.} (2006) have considered a general log-linear Poisson regression model to provide a comprehensive characterization of model choice and its uncertainty in time series studies of air pollution and mortality.  
 Enciso-Mora {\it et al.} (2009) have constructed a very efficient Reversible Jump Markov chain Monte Carlo (RJMCMC) algorithm to study the problem of model (order) selection for integer-valued autoregressive moving-average (INARMA) process. 
  More recently, Cleynen and Lebarbier (2014, 2017) have considered a particular case of model selection: multiple change-point detection. They considered segmentation problem for univariate distributions from the multiparameter exponential family, including integer-valued time series models. They have proposed a penalized log-likelihood estimator where the penalty function is constructed in a non-asymptotic context. See also Diop and Kengne (2019) for piecewise autoregression  in a general class of integer-valued time series.
 Alzahrani {\it et al.} (2018) have developed a procedure in a Bayesian framework for choosing between different classes of integer valued time series models.  
%
%

 \medskip

 In a model selection problem, the purpose is generally to construct a consistent or efficient procedure. 
 A consistent procedure leads to select the correct model (in a set of competing models, including the "true model"), with probability approaching one as the sample size increases to infinity. A procedure is said efficient when its risk is asymptotically equivalent to the risk of the oracle. 
 We address here the consistency question of model selection in integer-valued time series.  
 %
 We consider the general class of observation-driven models defined by
 
 \medskip
 {\bf Class} $\mathcal{OD}(f)$: A $\N_0$-valued ($\N_0=\N \cup \{0\}$) process $Y=\{Y_{t},~t\in \Z \}$ belongs to $\mathcal{OD}(f)$, if it satisfies:
   \begin{equation} \label{Model}
Y_t|\mathcal{F}_{t-1}\sim p(\cdot;\lambda_t) ~\text{ with } ~ \lambda_t=\E(Y_t|\mathcal{F}_{t-1})=f(Y_{t-1},Y_{t-2},\cdots),
\end{equation}
where $\mathcal{F}_{t-1}=\sigma\left\{Y_{t-1},Y_{t-2},\cdots \right\}$ (the $\sigma$-field generated by the whole past at time $t-1$),  
$p(\cdot)$ is a discrete distribution and $f(\cdot)$ is a measurable non-negative function defined on $\mathbb{N}_0^{\N}$.  

     \medskip
      \medskip 
We assume that the function $f(\cdot)$ is known up to a parameter $\theta \in \Theta$; where $\Theta$ is a compact set of $ \mathbb{R}^{d}$ ($d \in \N$). We consider then the class $\mathcal{OD}(f_\theta)$ given by

\medskip
{\bf Class} $\mathcal{OD}(f_\theta)$: A $\N_0$-valued process $Y=\{Y_{t},~t\in \Z \}$ belongs to $\mathcal{OD}(f_\theta)$, if it satisfies:
   \begin{equation} \label{Model_Bis}
Y_t|\mathcal{F}_{t-1}\sim p(\cdot;\lambda_t(\theta)) ~\text{ with }~ \lambda_t(\theta) = \mathbb{E}(Y_t|\mathcal{F}_{t-1})=f_\theta(Y_{t-1},\cdots).
\end{equation} 
 Each model of $ \mathcal{OD}(f_\theta)$ is associated to a distribution $p(\cdot)$ (assumed to be unknown) and a parameter $\theta \in \Theta$.
  Numerous classical time series models are included in $\mathcal{OD}(f_\theta)$: for instance Poisson, negative binomial, binary INGARCH or Poisson exponential autoregressive model (proposed by Fokianos {\it et al.} (2009)).  
This class of models has been studied by Ahmad and Francq (2016). 
 Under certain regularity conditions, they established the consistency and asymptotic normality of the Poisson quasi-maximum likelihood estimator of the model's parameter. 
 Cui and Qi (2017) have studied the inference in a particular case of model (\ref{Model_Bis}) where the conditional distribution belongs to the one-parameter exponential family.
 But, an asymptotic study of the model selection problem in the class $\mathcal{OD}(f_\theta)$ with unknown distribution $p(\cdot)$ and  infinite order processes (that provided a large way to take into account dependence on the past observations) has not yet been addressed.

\medskip
\medskip


%
%

\medskip
 Assume that the observations $Y_1,\cdots,Y_n$ are generated from the process $\{Y_{t},~t\in \Z\}$ satisfying (\ref{Model_Bis}), where the distribution $p(\cdot)$ and the true parameter (denoted by $\theta^*$) are unknown.
We consider a collection $\mathcal{M}$ of competing models belonging to $ \mathcal{OD}(f_\theta)$, containing the "true model" $m^*$, which corresponds to the parameter $\theta^*=\theta(m^*) \in \Theta$. 
Our main aim is to select the "best model" $\widehat{m}_n \in \mathcal{M}$ that we will consider as an estimator of $m^*$.  
%
 For all $m \in \mathcal{M}$, denote by $P_{n,\theta}$ the conditional distribution of $(Y_1,\cdots,Y_n)|\mathcal{F}_{n-1},m$  and consider the log-likelihood contrast given by:
\[
   ~ ~ \gamma_n( \theta|m ) := \gamma_n( P_{n,\theta} ) = - \log P_{n,\theta}(Y_1,\cdots,Y_n)  
             =  -  \sum_{t=1}^n \log p(Y_t; \lambda_t(\theta)),~~\forall~\theta \in \Theta(m),
  \]
 where $\Theta(m)$ is the parameter space of the model $m$, see below.
The estimator $\widehat{\theta}(m) $ of $\theta^*$ on the collection $\Theta(m)$ is obtained by minimizing the contrast $\gamma_n( \theta|m )$ over $\theta \in \Theta(m)$; that is
$ \widehat{\theta}(m) = \underset{ \theta \in \Theta(m)}{ \text{argmin}} ~\gamma_n( \theta|m )$.
%
  %
Thus, the estimator $\widehat m_n$ of $m^*$ is obtained by minimizing the penalized criterion 
\begin{equation}\label{Criterion-Pen}
\text{crit}_n(m) = \gamma_n( \widehat{\theta}(m)|m ) + \text{pen}_n(m), ~~\text{for all} ~ m \in \mathcal{M} ,
\end{equation}
where $\text{pen}_n : \mathcal{M} \rightarrow \R_+$ is a penalty function, possibly data-dependent. 
 Since the distribution $p(\cdot)$  is unknown, we propose a penalized criteria based on a Poisson quasi-likelihood and  provide sufficient conditions on the penalty $\text{pen}_n$, for which the estimator $\widehat{m}_n$ is consistent. 

\medskip
The paper is organised as follows. In Section 2, the assumptions and the definition of the Poisson quasi-maximum likelihood are provided. 
Section 3, we derive the model selection procedure and provide the main results. Some simulations results are displayed in Section 4, whereas Section 5 focus on applications on a real data example.  
 Section 6 is devoted to the proofs of the main results.
\section{Notation, Assumptions and Poisson QMLE}
\subsection{The framework}
 In the sequel, we will consider several models of $ \mathcal{OD}(f_\theta)$ and we define:
\begin{itemize}
	\item a model $m$ as a subset of $\{ 1,\cdots,d\}$ and denote by $|m|$ the dimension of model $m$;
	\item $\theta(m)$ as the parameter vector associated to the model $m$;
		\item $\Theta (m)=\{  (\theta_i)_{1 \leq i \leq d} \in \R^d,~~\theta_i=0 \text{ if } i \notin m \} \cap \Theta$ is  a compact set of $ \mathbb{R}^{d}$ ($d \in \N$) containing all of possible values for the parameter $\theta(m)$;
		\item $\mathcal{M}$ as a finite family of parametric models, $i.e$ $\mathcal{M} \in \mathcal{P}(\{ 1,\cdots,d\})$ (the set of all subsets of $\{ 1,\cdots,d\}$).
\end{itemize}
%
\begin{example}
 Assume that the observations $Y_1,\cdots,Y_n$ are generated from a negative binomial $INGARCH(p^*,q^*)$ process $\{Y_{t},~t\in \Z\}$; that is,
  \begin{equation}\label{NB_example}
      Y_{t}|\mathcal{F}_{t-1}\sim NB(r,p_t) \text{ with } r \frac{1 - p_t}{p_t} = \lambda_t= \alpha^*_{0} + \sum_{i=1}^{p^*} \alpha^*_i Y_{t-i} +  \sum_{i=1}^{q^*} \beta^*_i \lambda_{t-i} ;
  \end{equation}
  where $ \alpha_0^*>0$, $\alpha^*_1,\cdots, \alpha^*_{p^*} , \beta^*_1,\cdots,\beta^*_{q^*} \geq 0$, $\sum_{i=1}^{p^*} \alpha^*_i +  \sum_{i=1}^{q^*} \beta^*_i < 1$ and $NB(r,p)$ denotes the negative binomial distribution with parameters $r$ and $p$.
 The true parameter of the model is $\theta^* = (\alpha_0^*, \alpha^*_1,\cdots, \alpha^*_{p^*} , \beta^*_1,\cdots,\beta^*_{q^*})$.
 Since we can find  a sequence of non-negative real numbers
$(\psi_k(\theta^*))_{k \geq0}$ such that
\[ \lambda_t = \psi_0(\theta^*) + \sum_{k \geq1} \psi_k(\theta^*)Y_{t-k}, \]
 then, the model (\ref{NB_example}) belongs to the class $ \mathcal{OD}(f_{\theta^*})$.\\
 In practice, the conditional distribution of the observations is unknown (for this reason, we deal with the Poisson quasi-likelihood, see below); we will focus on the inference in the conditional mean $\lambda_t$. Based on the observations $Y_1,\cdots,Y_n$, the aim is to select the "best order" $(\widehat{p}_n,\widehat{q}_n)$ as an estimation of $(p^*,q^*)$.
 For this purpose, we consider a collection of $INGARCH(p,q)$ representation, with $(p,q) \in \{0,1,\cdots, p_{max} \} \times \{0,1,\cdots, q_{max} \}$ where $p_{max}, q_{max}$ are the fixed upper bounds of the orders, assumed to satisfy $p_{max} \geq p^*$ and $q_{max} \geq q^*$. 
 Therefore, consider $\Theta$ as a compact subset of $ (0,\infty) \times [0,\infty)^{p_{max} + q_{max}}$.
 Thus, a model $m$ is a subset of $\{ 1,\cdots,p_{max} + q_{max}\}$ with the parameter space 
 $\Theta (m)=\{  (\theta_i)_{1 \leq i \leq p_{max} + q_{max}+1} \in \R^{p_{max} + q_{max}+1},~~\theta_i=0 \text{ if } i \notin m \} \cap \Theta$.
\end{example}

\subsection{Notation and assumptions}
Throughout the sequel, the following norms will be used:
 {\em
\begin{itemize}
  \item $ \|x \|:= \sqrt{\sum_{i=1}^{p} |x_i|^2 } $, for any $x \in \mathbb{R}^{p}$; 


 \item  $\left\|f\right\|_{\Theta}:=\sup_{\theta \in \Theta}\left(\left\|f(\theta)\right\|\right)$ for any function $f:\Theta \longrightarrow \mathbb{R}^{d^{\prime}}$;

\item $\left\|Y\right\|_r:=\E\left(\left\|Y\right\|^r\right)^{1/r}$, if $Y$ is a random vector with finite $r-$order moments.
\end{itemize}
}

 \medskip
 
We set the following classical Lipschitz-type condition on the function $f_\theta$.

    \medskip
    \noindent \textbf{Assumption} \textbf{A}$_i (\Theta)$ ($i=0,1,2$):
    For any $y\in \mathbb{N}_{0}^{\N}$, the function $\theta \mapsto f_\theta(y)$ is $i$ times continuously differentiable on $\Theta$ 
    and
      there exists a sequence of non-negative real numbers $(\alpha^{(i)}_k)_{k\geq 1} $ satisfying
     $ \sum\limits_{k=1}^{\infty} \alpha^{(0)}_k <1 $ (or $ \sum\limits_{k=1}^{\infty} \alpha^{(i)}_k <\infty $ for $i=1, 2$);
   such that for any  $y, y' \in \mathbb{N}_{0}^{\N}$,
  \[ 
  \sup_{\theta \in \Theta  } \left| \frac{\partial^i f_\theta(y)}{ \partial\theta^i}-\frac{\partial^i f_\theta(y')}{ \partial\theta^i} \right|
  \leq  \sum\limits_{k=1}^{\infty}\alpha^{(i)}_k |y_k-y'_k|. 
  \]
%
    %
 
 \medskip 
 
 \noindent In the whole paper, it is assumed that there exists a stationary and ergodic process $\{Y_{t},~t\in \Z\}$ solution of (\ref{Model_Bis}), depending on the parameter $\theta^* \in \Theta(m^*)$ and satisfying
 \begin{equation}\label{moment}
    \exists C, \epsilon >0, \text{ such that } \forall t \in \Z, ~ ~  \E Y_{t}^{1+\epsilon} < C.
   \end{equation}

%
%

\subsection{Poisson QMLE}
We will carry out a brief overview of the Poisson quasi-maximum likelihood in the model (\ref{Model_Bis}) with the main asymptotic properties.  
 Consider that $(Y_1,\cdots,Y_n)$ is a trajectory generated from model (\ref{Model_Bis}) according to the true
parameter $\theta^* \in \Theta(m^*)$.  
 For any $m \in \mathcal{M}$ and $\theta \in \Theta(m)$, the Poisson quasi-log-likelihood is given by (up to a constant)

\[  L_n( \theta|m) := \sum_{t=1}^{n}(Y_t\log \lambda_t(\theta)- \lambda_t(\theta)) = \sum_{t=1}^{n} \ell_t(\theta) ~ \text{ with }
        \ell_t(\theta) = Y_t\log \lambda_t(\theta)- \lambda_t(\theta),
  \]
  where $ \lambda_t(\theta) =f_\theta(Y_{t-1},\cdots )$.
Since only $(Y_1, \cdots , Y_n)$ are observed, $L_n(\theta|m)$ is approximated by
\begin{equation}\label{logvm}
\widehat{L}_n(\theta|m) := \sum_{t=1}^{n}(Y_t\log \widehat{\lambda}_t(\theta)- \widehat{\lambda}_t(\theta)) = \sum_{t=1}^{n} \widehat{\ell}_t(\theta) ~ \text{ with }
 \widehat{\ell}_t(\theta) =  Y_t\log \widehat{\lambda}_t(\theta)- \widehat{\lambda}_t(\theta),
 \end{equation}
  where $ \widehat{\lambda}_t(\theta) =  \widehat{f}_\theta(Y_{t-1}, \cdots Y_{1},0,\cdots,0)$.\\
 The Poisson quasi-likelihood estimator (PQMLE) of $ \theta^*$ giving a model $m$ is defined by
 \begin{equation}\label{emv}
  \widehat{\theta}(m) :=  \underset{\theta\in \Theta(m)}{\text{argmax}} (\widehat{L}_n(\theta|m)).
  \end{equation}

 To ensure the consistency of the model selection procedure and the PQMLE, we impose the following regularity conditions on the "true model" $m^*$ (see also Ahmad and Francq \cite{Francq2016}):
 
\begin{enumerate}
   \item [(\textbf{A0}):] for all  $(\theta, \theta')\in \Theta^2$,
 $ \Big( f_\theta(y_{1}, y_{2}, \cdots)= f_{\theta'}(y_{1}, y_{2}, \cdots)  \ \text{a.s.} ~ \text{ for some } t \in \N \Big) \Rightarrow ~ \theta = \theta'$;
  moreover, $\exists  \underline{c}>0$ such that $\displaystyle \inf_{ \theta \in \Theta} f(y_{1}, y_{2}, \cdots; \theta)  \geq \underline{c}$ for all $ y \in  \N_0^{\N} $;
  \item [(\textbf{A1}):] $\theta^*$ is an interior point in the compact parameter space $\Theta(m^*) \subset \Theta \subset \mathbb{R}^{d}$;
  \item [(\textbf{A2}):] for all $y \in \mathbb{N}_{0}^{\N} $, the function $\theta \mapsto f_\theta(y)$ is twice continuously differentiable on $\Theta$; 
  \item [(\textbf{A3}):] $a_{t} \longrightarrow 0$ and $Y_{t} a_{t} \longrightarrow 0$ as $t\rightarrow \infty$, where $a_{t}=\underset{\theta \in \Theta }{\sup} \left| \widehat{\lambda}_{t}(\theta) -\lambda_{t}(\theta)\right|$;
   \item [(\textbf{A4}):] $J (\theta^*) =\E \Big[ \frac{1}{\lambda_{t}(\theta^* )}  \frac{\partial \lambda_{t}(\theta^* )}{ \partial \theta} \frac{\partial \lambda_{t}(\theta^* )}{ \partial \theta'}  \Big] <\infty$
           ~and~ $I(\theta^*) =\E \Big[ \frac{\var(Y_{t}|\mathcal{F}_{t-1})}{\lambda^2_{t}(\theta^* )}  \frac{\partial \lambda_{t}(\theta^* )}{ \partial \theta} \frac{\partial \lambda_{t}(\theta^* )}{ \partial \theta'}  \Big] <\infty$;
   \item [(\textbf{A5}):] for all $c' \in \R$, $c' \frac{\partial \lambda_{t} (\theta^*)}{\partial    \theta}=0$ a.s   $\Rightarrow ~ c'=0$ ;
   \item [(\textbf{A6}):] there exists a neighborhood $V(\theta^*)$ of $\theta^* $ such that: for all $i, j \in \left\{1,\cdots,d\right\} $, $$\E \left[ \sup_{\theta \in V(\theta^* )} \left|  \frac{\partial^2 }{ \partial \theta_i \partial \theta_j } \ell_{t}(\theta) \right|\right]<\infty ~;$$
   \item [(\textbf{A7}):]  $b_{t}$, $b_{t} Y_{t}$ and $a_{t}d_{t} Y_{t}$ are of order $O(t^{-h})$ for some $h>1/2$, where
    \[b_{t}=\underset{\theta \in \Theta }{\sup} \left\{\E \left[\left\|\frac{\partial \widehat{\lambda}_{t} (\theta)}{ \partial \theta} -\frac{\partial \lambda_{t}(\theta)}{ \partial \theta}  \right\| \right]\right\}~~\textrm{and}~~d_{t}=\underset{\theta \in \Theta }{\sup}\max \left\{ \E\left[\left\| \frac{1}{\widehat{\lambda}_{t}(\theta)} \frac{\partial \widehat{\lambda}_{t} (\theta)}{ \partial \theta} \right\|\right],  \E\left[ \left\| \frac{1}{\lambda_{t}(\theta)} \frac{\partial \lambda_{t} (\theta)}{ \partial \theta} \right\|\right] \right\}.  \]
    
\end{enumerate}
These assumptions hold for many classical linear and nonlinear models, we refer to Ahmad and Francq \cite{Francq2016}. 
\\
Under Assumptions (\textbf{A0})-(\textbf{A7}) and (\textbf{A}$_i (\Theta)$) (for $i=0,1,2$), Ahmad and Francq (2016) have established that the estimator $\widehat{\theta}(m^*)$ is consistent and asymptotically normal; that is,

\begin{equation}\label{Const_Theta}
\widehat{\theta}(m^*) \limitepsn \theta^*
\end{equation}
and
\begin{equation}\label{Normal_Theta}
\sqrt{n}\left(\widehat{\theta}(m^*)-\theta^*\right) \limiteloin \mathcal{N}_{|m^*|}(0,\Sigma_{\theta^*}) ,
\end{equation}
 where $\Sigma_{\theta^*}:=J^{-1}(\theta^*) I(\theta^*) J^{-1}(\theta^*)$.\\
%
%
 %
 %
  From the stationarity and ergodicity assumptions, for any  $m \in \mathcal{M}$ and $\theta \in \Theta(m)$, we have $\E \big[ \frac{1}{n} L_n(\theta | m) \big] = \E(\ell_1(\theta)).$ 
 Let us define 
  \begin{equation}\label{theta_etoile_m}
   \theta^*(m) :=  \underset{\theta\in \Theta(m)}{\text{argmax}} (\ell_1(\theta)). 
\end{equation} 
The following proposition establishes that, even in the misspecified framework, $\widehat{\theta}(m)$ is a consistent estimator of $\theta^*(m)$. 

 \begin{prop}\label{prop1}
 Assume that (\textbf{A0})-(\textbf{A3}), (\textbf{A}$_0(\Theta)$) and (\ref{moment}) (with $\epsilon > 1$) hold. Then:
\begin{enumerate}
\item For any model $m \in \mathcal{M}$,  there exists a unique $\theta^*(m) \in \Theta(m)$, satisfying (\ref{theta_etoile_m});
\item 
\begin{equation}\label{Const_Theta_m}
 \widehat{\theta}(m) \limitepsn \theta^*(m).
\end{equation}
\end{enumerate}

\end{prop}
 
  \section{ The model selection procedure and asymptotic results} 
In this section, we define the model selection criteria and provide the main 
 asymptotic results.  
 Assume that $\mathcal{M}$ is a finite family of the competing models satisfying the representation of class $\mathcal{OD}(f_\theta)$, including the "true model" $m^*$. 

\subsection{The penalized contrast}
%
For any model $m \in \mathcal{M}$, consider a contrast based on the Poisson quasi-likelihood: $-2 \widehat L_n(\widehat \theta(m))$. 
To select the "best model" $\widehat m_n$, we penalize this contrast by an additional term $\kappa_n |m|$, where $\kappa_n$ represents a regularization parameter. Define the penalized criteria by

\begin{equation}\label{Cont_pen}
 \widehat{C}_n(m):= -2 \widehat L_n(\widehat \theta(m)) + \kappa_n |m|, \text{ for all } m \in \mathcal{M} ,
 \end{equation}
 where
\begin{itemize}
	\item  $ (\kappa_n)_{n \in \N}$ is an increasing sequence satisfying $\kappa_n <n$ and $\kappa_n \rightarrow \infty$;
  \item $|m|$ is the number of estimated components of $\theta \in \Theta(m)$, also called the dimension of the model $m$. 
\end{itemize}
	The choice of the "best model" $\widehat m_n$ is done by minimizing the penalized criteria:
\begin{equation}\label{Estim_m}
  \widehat m_n :=  \underset{  m \in \mathcal{M}}{\text{argmin}} \left(\widehat{C}_n(m)\right).
  \end{equation}
  
  \subsection{Asymptotic results}
In the sequel, the competing collection of models $\mathcal{M}$ is assumed to be finite.
Let us consider $\widehat m_n$, a selected model with respect to $\mathcal{M}$. 
The following theorem establishes the consistency of the model selection procedure.
\begin{thm}\label{Th1}
Assume that (\textbf{A0})-(\textbf{A7}),  (\textbf{A}$_i(\Theta)$) (for $i=0,1,2$) and (\ref{moment}) (with $\epsilon > 1$) hold. 
If $(\kappa_{n})$ satisfies
\begin{equation}\label{cond_th1}
 \sum_{\ell \geq 1} \frac{1}{ \kappa_{\ell} } \sum_{k \geq \ell } \alpha_k^{(i)}   < \infty, \text{ for } i=0,1,2
 \end{equation}
 and $\kappa_n=o(n)$, then
\begin{equation}\label{Cons_C}
\mathbb{P}(\widehat m_n =m^*) \limiten 1 ~~and ~~   \widehat \theta(\widehat m_n) \limiteproban \theta^*.
\end{equation}
\end{thm}
The next theorem shows the asymptotic normality of the PQMLE of the chosen model.
\begin{thm}\label{Th2} 
Assume that (\textbf{A0})-(\textbf{A7}),  (\textbf{A}$_i(\Theta)$) (for $i=0,1,2$), (\ref{moment}) (with $\epsilon > 1$) and (\ref{cond_th1}) hold. 
 If $\kappa_n=o(n)$, then
\begin{equation}\label{Normal_C}
 \sqrt{n}\left(\widehat{\theta}(\widehat{m}_n)-\theta^*\right) \limiteloin \mathcal{N}_d \left(0,\Sigma \right), 
\end{equation}
 where $\Sigma:=J^{-1}(\theta^*) I (\theta^*) J^{-1}(\theta^*)$ with
\[
J (\theta^*)=\E \Big[ \frac{1}{\lambda_{1}(\theta^* )}  \frac{\partial \lambda_{1}(\theta^* )}{ \partial \theta} \frac{\partial \lambda_{1}(\theta^* )}{ \partial \theta'}  \Big]
           ~~ \textrm{and}~~  I (\theta^* )=\E \Big[ \frac{\var(Y_{1}|\mathcal{F}_{0})}{\lambda^2_{1}(\theta^* )}  \frac{\partial \lambda_{1}(\theta^* )}{ \partial \theta} \frac{\partial \lambda_{1}(\theta^* )}{ \partial \theta'}  \Big].
 \]
\end{thm}

 \medskip
  \medskip
\noindent Let us stress that the condition (\ref{cond_th1}) on $(\kappa_n)_{n\in \N}$ can be easily obtained if the Lipschitzian coefficients of $f_\theta(\cdot)$ and their derivatives are bounded by a geometric or Riemanian decrease:
 %
\begin{enumerate} 
	\item  the geometric case: if $\alpha^{(i)}_k= O(a^k)$ ($i=0,1,2$) with $0 \leq a <1$, then any choice of $(\kappa_n)_{n\in \N}$ such that $\kappa_n =o(n) $ and $\kappa_n \rightarrow \infty$ satisfies (\ref{cond_th1}) ; and the consistency holds (for instance, the BIC approach given by $\kappa_n=\log n$).
 \item the Riemanian case: if $\alpha^{(i)}_k= O(k^{-\gamma})$ ($i=0,1,2$) with $\gamma >0$,
\begin{itemize}
	\item if $\gamma >2$, then the condition (\ref{cond_th1}) holds for any choice of $(\kappa_n)_{n\in \N}$ such that $\kappa_n  =o(n)$ and $\kappa_n \rightarrow \infty$.
	\item if $0<\gamma \leq 2$, then one can choose any sequence such that $\kappa_n=O(n^\delta)$ with $\delta >2-\gamma$.
\end{itemize}
\end{enumerate}


\section{Some simulations results}

This section presents a simulation study to illustrate the performances of the selection procedure proposed.    
 We will compare the performances of the penalties $\log n$ and $n^{1/3}$ for a linear and nonlinear dynamic models. 
  For each model considered, we use a Monte Carlo experiments with the sample size $n$ belongs to $\{ 500,\,1000,\, 2000 \}$.

\subsection{Linear dynamic models}

We consider the following models:
\begin{enumerate}
 	\item \emph {Poisson-INARCH$(2)$ process}: 
	\begin{equation*}
                   \text{Model \textsl{A.} : }  	Y_{t}|\mathcal{F}_{t-1} \sim \mathcal{P}(\lambda_{t})~~\text{with}~~  \lambda_{t}= 
                                                	0.5+ 0.3 Y_{t-1}+ 0.25 Y_{t-2} ;
    \end{equation*}
    \item \emph {Poisson-INARCH$(1,1)$ process}: 
    \begin{equation*}
	                   \text{Model \textsl{B.} : } Y_{t}|\mathcal{F}_{t-1} \sim \mathcal{P}(\lambda_{t})~~\text{with}~~                                         \lambda_{t}=  1 + 0.3 Y_{t-1}+ 0.45\lambda_{t-1} ;
    \end{equation*}
    \item \emph {Binary-INARCH$(2)$ process}: 
    \begin{equation*}
	             \text{Model \textsl{C.} : } Y_{t}|\mathcal{F}_{t-1} \sim \mathcal{B}(p_{t})~~\text{with}~~  p_{t}= \lambda_{t} =
	                                       0.15 + 0.25 Y_{t-1}+ 0.2 Y_{t-2} ;
    \end{equation*}
   \item \emph {Binary-INGARCH$(1,1)$ process}: 
    \begin{equation*}
                	      \text{Model \textsl{D.} : }  Y_{t}|\mathcal{F}_{t-1} \sim \mathcal{B}(p_{t})~~\text{with}~~  p_{t}= \lambda_{t} =
                                       	0.1 + 0.35 Y_{t-1}+ 0.4 \lambda_{t-1} ;
    \end{equation*} 
\end{enumerate}
 where $\mathcal{P}(\lambda)$ is the Poisson distribution with parameter $\lambda$ and $\mathcal{B}(p)$ the Bernoulli distribution with parameter $p$. 
We consider as competitive models all INGARCH$(p,q)$ representations with $ (p,q) \in \{0,\cdots,  5\} \times \{0,\cdots,  5\}$. The results of the selection procedure are presented in Table \ref{Res_sim1}.  

 \begin{table}[h!]
\footnotesize
\centering
\caption{\it Frequencies of selected model based on $100$ replications depending on sample's length for Model \textsl{A.},  \textsl{B.}, \textsl{C.} and  \textsl{D.}}
\label{Res_sim1}
\vspace{.2cm}
\begin{tabular}{clcccccccc}

\Xhline{.9pt}
 & &&&& \\
& &  \multicolumn{2} {c} {$n=500$}&&\multicolumn{2} {c} {$n=1000$}&&\multicolumn{2} {c} {$n=2000$}   \\
\cline{3-4}\cline{6-7}\cline{9-10}
 &&$\log n$&$n^{1/3}$ &&$\log n$&$n^{1/3}$ &&$\log n$&$n^{1/3}$ \\
\Xhline{.72pt}
                  & $|\widehat m_n|<|m^*|$      &$0.00$&$0.00$  &&$0.00$&$0.00$  &&$0.00$&$0.00$  \\
Model \textsl{A.} & $\widehat m_n=m^*$          &$0.92$&$0.92$  &&$0.95$&$0.95$  &&$1.00$&$1.00$  \\
                  & $|\widehat m_n| \geq |m^*|$ and  $\widehat m_n \neq m^*$ &$0.08$&$0.08$  &&$0.05$&$0.05$  &&$0.00$&$0.00$  \\
                      

 \hline
 \hline
                   &  $|\widehat m_n|<|m^*|$       &$0.12$&$0.19$  &&$0.00$&$0.00$  &&$0.00$&$0.00$ \\
 Model \textsl{B.} & $\widehat m_n=m^*$            &$0.70$&$0.67$  &&$0.95$&$0.95$  &&$0.96$&$0.97$ \\
                   & $|\widehat m_n| \geq |m^*|$ and  $\widehat m_n \neq m^*$ &$0.18$&$0.14$  &&$0.05$&$0.05$  &&$0.04$&$0.03$  \\
          
 \hline
 \hline
   
                   & $|\widehat m_n|<|m^*|$       &$0.17$&$0.24 $  &&$0.00$&$0.01$  &&$0.00$&$0.00$  \\
Model \textsl{C.}  & $\widehat m_n=m^*$           &$0.72$&$0.66$  &&$0.91$&$0.90$  &&$0.98$&$0.98$  \\
                   & $|\widehat m_n| \geq |m^*|$ and  $\widehat m_n \neq m^*$ &$0.11$&$0.10$  &&$0.09$&$0.09$  &&$0.02$&$0.02$  \\
           
\hline
 \hline   
                   & $|\widehat m_n|<|m^*|$       &$0.38$&$0.44$  &&$0.00$&$0.03$  &&$0.00$&$0.01$  \\
Model \textsl{D.}  & $\widehat m_n=m^*$           &$0.49$&$ 0.43$  &&$0.87$&$0.85$  &&$0.95$&$0.94$  \\
                   & $|\widehat m_n| \geq |m^*|$ and  $\widehat m_n \neq m^*$ &$0.13$&$0.13$  &&$0.13$&$0.12$  &&$0.05$&$0.05$  \\

 \Xhline{.9pt}
\end{tabular}
\end{table}
 The results of Table \ref{Res_sim1} show that for both the penalties, the performances of the procedure increase with $n$ for all models. We can see that the consistency of the penalties $\log n$ and $n^{1/3}$ is numerically convincing, which is in accordance with the asymptotic results of Theorem \ref{Th1}. 
 However, the performances are more interesting for the models without moving average component. 
    We also note that for a small sample size (see for instance, $n=500$), the $\log n$-penalty  is a little bit better in comparison with the $n^{1/3}$-penalty except for the Poisson-INARCH$(2)$ process.

\subsection{Nonlinear dynamic models}

Consider the autoregressive model defined by 

\begin{equation}\label{NB_Nonlinear_model}
	                     Y_{t}|\mathcal{F}_{t-1} \sim NB(r,p_{t})~~\text{with}~~  r\frac{(1-p_{t})}{p_{t}}= \lambda_{t} =
	                                       f_{\theta^*}(Y_{t-1},\lambda_{t-1});
    \end{equation}
   where $r \in \N$, $\theta^*$ is the true parameter belonging in a compact set $\Theta$, for any $\theta \in \Theta $, $f_\theta$ is a non-negative nonlinear measurable  function defined on $\N_0 \times \R_{+}$ and $NB(r, p)$ represents the negative binomial distribution with parameters $r$ and $p$. We assume that the function $f_\theta$ satisfies the contraction condition;  {\it i.e.},   there exist non-negative constants $a$ and $b$ satisfying $a + b < 1$ such that for any $(y, \lambda) \in \N_0 \times \R_{+}$ and $(y', \lambda') \in \N_0 \times \R_{+}$,
\begin{equation}\label{nonlinear_lip}
\left\|f_\theta(y,\lambda)- f_\theta(y',\lambda')\right\|_{\Theta}  \leq a \left|y-y' \right| + b\left|\lambda-\lambda' \right|.
\end{equation}
\indent According to (\ref{nonlinear_lip}), for any $\theta \in \Theta$, we can find a measurable function
 $f^{\theta^*}_\infty:[0,\infty)^{\N}\rightarrow [0,\infty) $  such that
 \[  \lambda_t=f^{\theta^*}_\infty(Y_{t-1},Y_{t-2},\cdots); \]
 that is, the model (\ref{NB_Nonlinear_model}) belongs to the class $\mathcal{OD}(f_{\theta^*})$.
 Moreover, the process $\{Y_t,~t\geq1 \}$ is absolutely regular with geometrically decaying coefficients and $\{(Y_t,\lambda_t),~t\geq1 \}$ is strictly stationary and ergodic (see \cite{Davis2016}). In case where the function $f_\theta$ is linear   (with INGARCH(1,1) representation as in (\ref{NB_example})), the existence of the second-order moment of $Y_t$ has already been studied (see for instance, Christou and Fokianos (2014)). 
The following proposition establishes the existence of the second-order moment of $Y_t$ in case where $f_\theta$ is nonlinear (model (\ref{NB_Nonlinear_model})).

\begin{prop}\label{prop_non_lin} 
Assume that (\ref{nonlinear_lip}) holds. A sufficient condition for that $\E(Y^2_t)<\infty$ for all $t \in \Z$ is 
\begin{equation}\label{cond_prop_nonl}
(a+b)^2 +\frac{a^2}{r} <1.
\end{equation}
\end{prop}  
 \medskip 
\medskip 
Now, consider the particular case of model (\ref{NB_Nonlinear_model}) given by
\begin{equation}\label{model_noeud}
    Y_{t}|\mathcal{F}_{t-1} \sim NB(r,p_{t})~~\text{with}~~  r\frac{(1-p_{t})}{p_{t}}= \lambda_{t} =
    \alpha^*_{0} + \alpha^*_{1}Y_{t-1} + \alpha^*_{2}\lambda_{t-1} +  \sum_{k =1}^{K^*} \beta^*_{k}(Y_{t-1}-\xi^*_k)^+ , 
   \end{equation}
 where $K^* \in\N_0$, $\alpha^*_{0}>0$, $\alpha^*_{i} \geq 0$ (for $i=1,2$), $\beta^*_{k} \geq 0$ (for $k=1,\cdots,K$),   
 $~(\xi^*_k)_{1 \leq k \leq K}$ are non-negative integers (so-called knots) and $x^+ = \max(x,0)$ is the positive part of $x$.  
 This process is a special case of (\ref{NB_Nonlinear_model}), where $\theta^*=(\alpha^*_{0},\alpha^*_{1},\alpha^*_{2},\beta^*_{1},\cdots,\beta^*_{K^*})$  and $f_{\theta^*}(y,\lambda)= \alpha^*_{0} + \alpha^*_{1}y +\alpha^*_{2}\lambda + \sum_{k =1}^{K^*} \beta^*_{k}(y-\xi^*_k)^+$, for any $(y,\lambda) \in \N_{0} \times \R_{+}$.     
 In particular, when $K^*=0$, model (\ref{model_noeud}) reduces to the NB-INGARCH$(1,1)$. The inference question in this model have also been studied by Davis and Liu (2012).  
The aim is to select the "best" number of knots in $\{0,1,\cdots, K_{max}\}$; where  $K_{max} \in \N_0$ is a fixed upper bound, assumed to satisfy $K_{max} \geq K^*$. 
Thus, $\Theta$ is a compact subset of $(0,\infty) \times [0,\infty)^{K_{max} + 2}$ such that $\alpha_{1}+\alpha_{2}+\sum_{k =1}^{K} \beta_{k}<1$ for all $\theta = (\alpha_{0},\alpha_{1},\alpha_{2},\beta_{1},\cdots,\beta_{K}) \in \Theta$. The true parameter $\theta^*$ could be rewritten as $\theta^*=(\alpha^*_{0},\alpha^*_{1},\alpha^*_{2},\beta^*_{1},\cdots,\beta^*_{K^*},0,\cdots,0)$, so it is an element of $\Theta$. 
%
%
 %
 
  \medskip
 For $r=1$ and $r=8$, we generate a trajectory of model (\ref{model_noeud}) with the following parameters: $K^*=1$, $\xi^*_1=2$ and $\theta^*=(1,0.2,0.15,0.35)$. The competing models considered are all process satisfying (\ref{model_noeud}) with $ K \in \{0,1,2,3\}$ and $\xi_k \in \{1,2,3,4\}$ (for any $k=1,\cdots,K$). We will focus on the selection of the "best value" of $K^*$ (denoted by $\widehat K_n$), which allows to determine the dimension of the model.  
 Table \ref{Res_sim2} indicates the frequencies of number of replications where $\widehat K_n<K^*$, $\widehat K_n=K^*$ and $\widehat K_n>K^*$. 

 \begin{table}[h!]
\footnotesize
\centering
\caption{\it Frequencies of the selection of the true, low and high value of $K^*$ for model (\ref{model_noeud}) based on $100$ replications.}
\label{Res_sim2}
\vspace{.2cm}
\begin{tabular}{cclccccccc}

\Xhline{.9pt}
  &&&&&& \\
 &&  \multicolumn{2} {c} {$n=500$}&&\multicolumn{2} {c} {$n=1000$}&&\multicolumn{2} {c} {$n=2000$}   \\
\cline{3-4}\cline{6-7}\cline{9-10}
  &&$\log n$&$n^{1/3}$ &&$\log n$&$n^{1/3}$&&$\log n$&$n^{1/3}$ \\
\Xhline{.72pt}
       &$\widehat K_n<K^*$  &$0.36$&$0.42$  &&$0.13$&$0.27$  &&$0.02$&$0.07$  \\
 $r=1$ &$\widehat K_n=K^*$  &$0.64$&$0.58$  &&$0.87$&$0.73$  &&$0.94$&$0.93$  \\
       &$\widehat K_n>K^*$  &$0.00$&$0.00$  &&$0.00$&$0.00$  &&$0.04$&$0.00$ \\                     

 \hline
 \hline
        &$\widehat K_n<K^*$  &$0.35$&$0.48$  &&$0.12$&$0.21$  &&$0.01$&$0.04$  \\
 $r=8$ &$\widehat K_n=K^*$   &$0.65$&$0.52$  &&$0.88$&$0.79$  &&$0.99$&$0.96$  \\
       &$\widehat K_n>K^*$   &$0.00$&$0.00$  &&$0.00$&$0.00$  &&$0.00$&$0.00$ \\ 
       
       \Xhline{.9pt}
\end{tabular}
\end{table}

Once again, Table \ref{Res_sim2} shows that the performances of the proposed procedure increase with the sample size. Although the $\log n$-penalty outperforms the $n^{1/3}$-penalty for moderate sample size (when $n=500$ and $n=1000$), the performances displayed by theses penalties are close when $n=2000$. Moreover, the empirical evidence of the consistency of the proposed procedure appears to be quite convincing. 
 
%
%
 
 \section{Real data application} 
\begin{figure}[h!]
\begin{center}
\includegraphics[height=6cm, width=11cm]{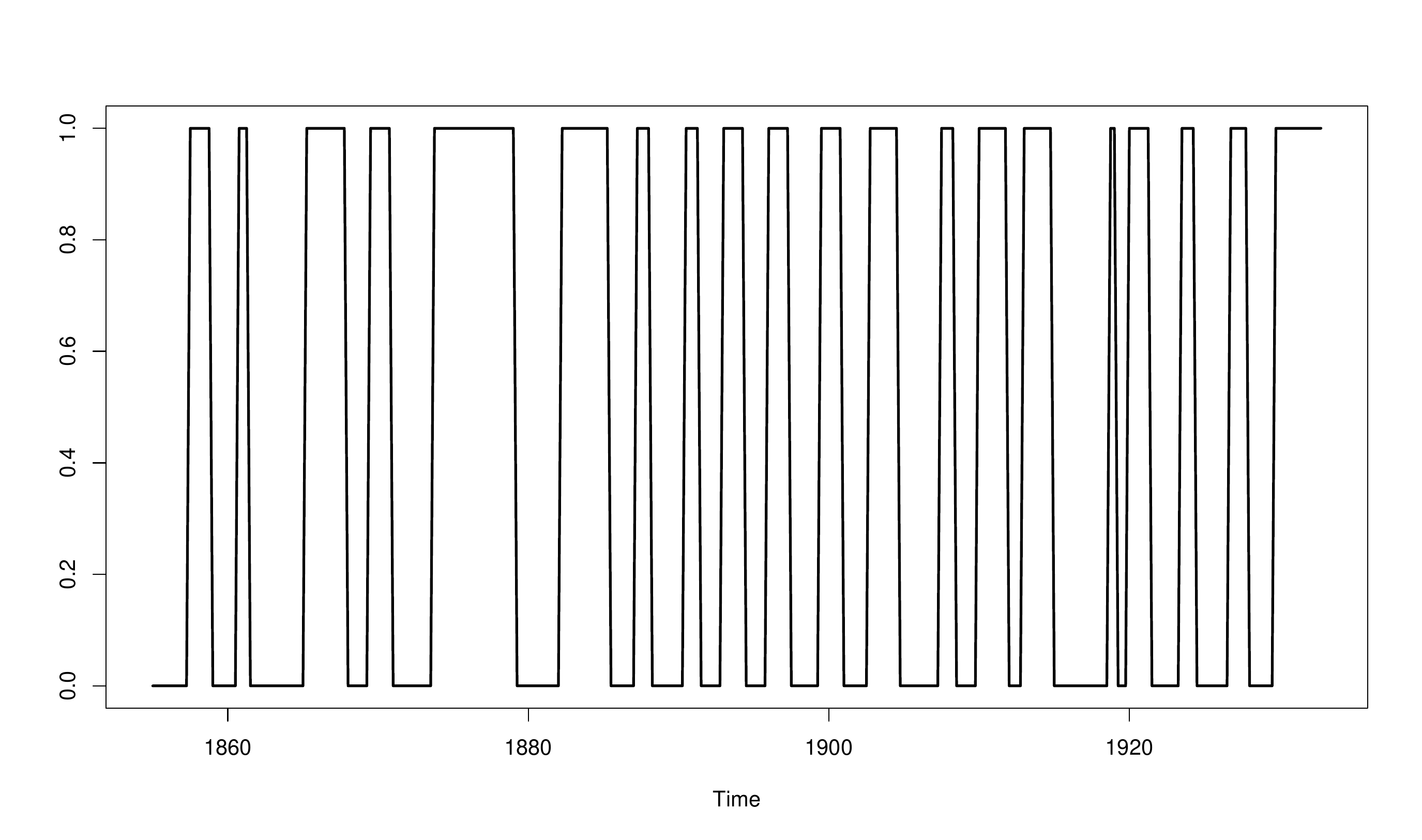}
\end{center}
\vspace{-.8cm}
\caption{ The USA recession data in the period 1855-1932.}
\label{Graphe_USA}
\end{figure}

 We apply the proposed procedure to the quarterly recession data from the USA. There are $636$ observations from $1855$ to $2013$, available at "https://fred.stlouisfed.org/series/USREC". The series $(Y_t)$ is a binary variable that is equal to $1$ if there is a recession in at least one month in the quarter and $0$ otherwise. In the literature, several works have already been carried out on these data (see for instance, Startz  (2008)). 
 In the context of break detection, Diop and Kengne (2017) have analyzed these data by applying a change-point test based on the maximum likelihood estimator of the model's parameter. They have detected a break (two regimes) in the last quarter of 1932 (at the time $t=312$).  
 Here, we limit ourselves to the first regime; \textit{i.e.}, the first $312$ observations which are represented in Figure \ref{Graphe_USA}.
 
 \medskip
 
 We now consider the collection of all INGARCH$(p,q)$ representation with $(p,q) \in \{0,\cdots, 5\} \times \{0,\cdots, 5\}$, which leads us to $36$ competing models. To select the "best" orders $\widehat{p}$ and $\widehat{q}$, we apply the selection procedure based on the Poisson quasi-likelihood with penalties $\log n$ and $n^{1/3}$. The obtained results show that the INARCH$(1)$ (\textit{i.e.}, $\widehat{p}=1$ and $\widehat{q}=0$) representation is the "best" model according to both criteria. This is in accordance with the conclusions of Diop and Kengne (2017), where their procedure is based on the maximum likelihood estimator. The estimated model with the PQMLE is  
 \begin{equation}\label{estim_USA}
 \begin{array}{l}
 \E(Y_t|\mathcal{F}_{t-1})=0.120 + 0.748 Y_{t-1},\\
       ~~~~~~~~~~~~~~~~~ (0.029)~~(0.216) 
 \end{array}
 \end{equation}
 where in parentheses are the standard errors of the estimators obtained from the robust sandwich matrix. \\ 
 Let us stress that for the sample sizes $n<500$, the results have not been presented in the simulation experiments. But, for $n=300$, we have carried out the binary-INARCH(1) model with the same scenario (as in (\ref{estim_USA})) and the numerical results show that the frequency of choosing the true model is very close to $100\%$. This lends a substantial support to the representation retained by the selection procedure to fit these data.
 
 \section{Proofs of the main results}

 In the sequel, we set $f^\theta_{t}:=f_\theta(Y_{t-1},\cdots)$. Also, $C$ denotes a positive constant whom value may differ from an inequality to another.
 
 \subsection{Proof of Proposition \ref{prop1}}
 \begin{enumerate}
 \item Let $m \in \mathcal{M}$. Consider the  function $L^m : \Theta(m)  \rightarrow \R$, defined by $L^{m}(\theta) =\E\left[\ell_{1}(\theta)\right]$, for all $\theta \in \Theta(m)$. 
For any $\theta \in \Theta(m) $, we have
\begin{align*}
L^{m}(\theta^*(m)) - L^{m}(\theta) &= \E[\ell_{1}(\theta^*(m))] -\E[\ell_{1}(\theta)] \\
&= \E[\ell_{1}(\theta)] -\E[\ell_{1}(\theta^*(m))]\\
&= \E[Y_{1}\log f^{\theta ^*(m)}_{1}- f^{\theta ^*(m)}_{1}] -\E[Y_{1}\log f^\theta_{1}- f^\theta_{1}]\\
&= \E \Big[ f^{\theta ^*(m)}_{1} \big(\log f^{\theta^*(m)}_{1} -\log f^\theta_{1} \big) \Big] - \E \Big[f^{\theta ^*(m)}_{1} - f^\theta_{1} \Big].
\end{align*}
By applying the mean value theorem to the function $x \mapsto \log x$ defined in $[\underline{c},+\infty[$;
there exists $\xi$ between $f^{\theta ^*(m)}_{1}$ and $f^{\theta}_{1}$ such that
 \[
 \log f^{\theta ^*(m)}_{1} -\log f^\theta_{1}= \frac{1}{\xi}\left(f^{\theta ^*(m)}_{1} - f^\theta_{1}\right).
 \]
 Hence,
 \begin{align*}
L^{m}(\theta^*(m)) - L^{m}(\theta)
&= \E \Big[ \frac{f^{\theta ^*(m)}_{1}}{\xi}\left(f^{\theta ^*(m)}_{1} - f^\theta_{1}\right) \Big] - \E \Big[f^{\theta ^*(m)}_{1} - f^\theta_{1} \Big]\\
&= \E \Big[ \Big( \frac{f^{\theta ^*(m)}_{1}}{\xi} -1 \Big)\left(f^{\theta ^*(m)}_{1} - f^\theta_{1}\right) \Big]\\
&= \E \Big[ \frac{1}{\xi}\left( f^{\theta ^*(m)}_{1}-\xi  \right)\left(f^{\theta ^*(m)}_{1} - f^\theta_{1}\right) \Big].
\end{align*}
From Assumption (\textbf{A0}), it follows that  $\frac{1}{\xi}\left( f^{\theta ^*(m)}_{1}-\xi  \right)\left(f^{\theta ^*(m)}_{1} - f^\theta_{1}\right) \neq 0 $ a.s.,  if  $\theta \neq \theta ^*(m)$.\\
Moreover,
\begin{itemize}
    \item if $f^{\theta ^*(m)}_{1} <f^{\theta}_{1}$, then $f^{\theta ^*(m)}_{1}< \xi <f^{\theta}_{1}$ and hence $\frac{1}{\xi}\left( f^{\theta ^*(m)}_{1}-\xi  \right)\left(f^{\theta ^*(m)}_{1} - f^\theta_{1}\right) > 0 $;

    \item if $f^{\theta ^*(m)}_{1} >f^{\theta}_{1}$, then $f^{\theta}_{1}< \xi <f^{\theta^*(m)}_{1}$ and hence $\frac{1}{\xi}\left( f^{\theta ^*(m)}_{1}-\xi  \right)\left(f^{\theta ^*(m)}_{1} - f^\theta_{1}\right) > 0 $.
\end{itemize}
We deduce that  $\frac{1}{\xi}\left( f^{\theta ^*(m)}_{1}-\xi  \right)\left(f^{\theta ^*(m)}_{1} - f^\theta_{1}\right) > 0 $ a.s.. Hence, $L^{m}(\theta^*(m)) - L^{m}(\theta)>0$ a.s., if $\theta \neq \theta^*(m)$. 
Thus, the  function $L^{m}(\theta)$ has a unique maximum at $\theta^*(m)$.
\item Recall that, since $\{Y_{t},~t\in \Z\}$ is stationary and ergodic, for any $m \in \mathcal{M}$ and $\theta \in \Theta(m)$, the process $\{\ell_{t}(\theta),~t\in \Z\}$ is also a stationary and ergodic sequence.
 Let us show that $\left\| \frac{1}{n}\widehat{L}_{n}(\theta)-\E (\ell_{1}(\theta))\right\|_{\Theta(m)} \limiten 0$.
 %
  Recall that $\ell_{t}(\theta) = Y_{t}\log f^\theta_{t}- f^\theta_{t}$, for any $\theta \in \Theta(m)$. We have
\begin{align*}
\left|\ell_{t} \right| &\leq Y_{t} \left|\log f^\theta_{t}\right| + \left|f^\theta_{t}\right|\\
&\leq
 Y_{t} \Big|\log \Big(\frac{f^\theta_{t}}{\underline{c}} \times \underline{c}\Big)\Big| + \left|f^\theta_{t}\right|\\
&\leq
 Y_{t} \left(\Big|\frac{f^\theta_{t}}{\underline{c}} -1\Big|  +\left| \log \underline{c} \right|\right) + \left|f^\theta_{t}\right|~~(\text{because for}~x>1,~\left|\log x\right| \leq \left| x-1\right|)\\
 &\leq
 Y_{t} \left(\Big|\frac{f^\theta_{t}}{\underline{c}}\Big| +1  +\left| \log \underline{c} \right|\right) + \left|f^\theta_{t}\right|.
\end{align*}
Hence,
\begin{equation}\label{Major_N_ell_m}
\left\| \ell_{t} \right\|_{\Theta(m)} \leq
 Y_{t} \left(\frac{1}{\underline{c}} \left\|f^\theta_{t}\right\|_{\Theta(m)} +1  +\left| \log \underline{c} \right|\right) + \left\|f^\theta_{t}\right\|_{\Theta(m)} .
\end{equation}
We will show that $\E \left[\left\| \ell_{t} \right\|_{\Theta(m)}\right] < \infty$. According to (\ref{Major_N_ell_m}), we have
\begin{align*}
\E \left[\left\| \ell_{t} \right\|_{\Theta(m)}\right] & \leq
\E \left[Y_{t} \left(\frac{1}{\underline{c}} \left\|f^\theta_{t}\right\|_{\Theta(m)} +1  +\left| \log \underline{c} \right|\right) + \left\|f^\theta_{t}\right\|_{\Theta(m)}  \right]\\
&\leq
C \E \left[\left(\frac{Y_{t}}{\underline{c}} +1\right) \left\|f^\theta_{t}\right\|_{\Theta(m)}  \right]\\
&\leq
  C\left(\E\left[ \left(\frac{Y_{t}}{\underline{c}}  +1\right)^2 \right]\right)^{1/2} \times \left(\E \left\|f^\theta_{t} \right\|^2_{\Theta(m)} \right)^{1/2} \\
&\leq
C \left(\E \left\|f^\theta_{t} \right\|^2_{\Theta(m)}\right)^{1/2}.
\end{align*}
In addition, according to (\textbf{A}$_0 (\Theta)$), we have
\begin{align*}
\left\|f^\theta_{t}\right\|_{\Theta(m)} & \leq \left\|f^\theta_{t} - f^\theta(0,\cdots) \right\|_\Theta + \left\|f^\theta(0,\cdots) \right\|_{\Theta(m)}
 \leq  \sum\limits_{k \geq 1} \alpha^{(0)}_k |Y_{t-k}|  +\left\|f^\theta(0,\cdots) \right\|_{\Theta(m)}.
\end{align*}
Therefore,
\begin{align*}
\E \left[\left\| \ell_{t} \right\|_{\Theta(m)}\right] & \leq
C \Big[\E \Big( \sum\limits_{k \geq 1} \alpha^{(0)}_k |Y_{t-k}|  +\left\|f^\theta(0,\cdots) \right\|_\Theta \Big)^2_\Theta\Big]^{1/2}\\
&\leq
C   \sum\limits_{k \geq 1} \alpha^{(0)}_k \left(\E |Y_{t-k}|^2\right)^{1/2}  +\left(\E \left\|f^\theta(0,\cdots) \right\|^2_\Theta \right)^{1/2}\\
&\leq
C   \sum\limits_{k \geq 1} \alpha^{(0)}_k +\left(\E \left\|f^\theta(0,\cdots) \right\|^2_{\Theta(m)} \right)^{1/2} <\infty.
\end{align*}

By the uniform strong law of large number applied on the process  $\{\ell_{t}(\theta),~t\in \Z\}$, it holds that
\begin{equation} \label{conv_Ln_m}
\left\| \frac{1}{n} L_{n}(\theta)-\E (\ell_{1}(\theta))\right\|_{\Theta(m)}= \left\| \frac{1}{n} \sum_{t=1}^{n} \ell_{t}(\theta)-\E (\ell_{1}(\theta))\right\|_{\Theta(m)}  \limiten 0.
 \end{equation}
Moreover, by going along similar lines as \cite{Francq2016} (see also the proof of Lemma \ref{lem0} below), we get 
\begin{equation}\label{conv1_Francq}
 \frac{1}{n}\left\| \widehat{L}_{n}( \theta)-L_{n}(\theta)\right\|_{\Theta(m)} \limiten 0.
 \end{equation}
From (\ref{conv_Ln_m}) and (\ref{conv1_Francq}), we deduce that
\begin{equation} \label{conv_Ln_chapo_m}
\left\| \frac{1}{n} \widehat{L}_{n}(\theta)-\E (\ell_{1}(\theta))\right\|_{\Theta(m)} \limiten 0.
 \end{equation}
\noindent Thus, the previous item, (\ref{conv_Ln_chapo_m}) and standard arguments lead to the consistency of $\widehat{\theta}(m)$.  
 \end{enumerate}
$~~~~~~~~~~~~~~~~~~~~~~~~~~~~~~~~~~~~~~~~~~~~~~~~~~~~~~~~~~~~~~~~~~~~~~~~~~~~~~~~~~~~~~~~~~~~~~~~~~~~~~~~~~~~~~~~~~~~~~~~~~~~~~~~~~~~~~~~~~~~~~~~~~~~ \blacksquare$

 \subsection{Proof of Theorem \ref{Th1}}
 
 Throughout this section, we consider the following lemma. 
 
 \begin{lem}\label{lem0} 
 Suppose that the assumptions of Theorem \ref{Th1} hold. Then
 
\begin{equation}\label{eq_lem0}
\frac{1}{\kappa_n} \left\|\widehat L_n(\theta) - L_n(\theta)  \right\|_\Theta \limitepsn 0 .
 \end{equation}
 \end{lem}
 \emph{\bf Proof of Lemma \ref{lem0}}\\
 Remark that, for all $n \in \N$,
\begin{align*}
  \| \widehat{L}_n(\theta)-L_n(\theta)\|_\Theta
 &\leq
 \sum_{t=1}^{n}\|\widehat{\ell}_{t}(\theta)-\ell_t(\theta) \|_\Theta\\
 &\leq
  \sum_{t=1}^{n} \|Y_t\log \widehat{\lambda}_{t}(\theta)- \widehat{\lambda}_{t}(\theta) -Y_t\log \lambda_t(\theta)+ \lambda_t(\theta)\|_\Theta\\
  &\leq
  \sum_{t=1}^{n}  (Y_t \|\log \widehat{f}^\theta_t-\log f^\theta_t\|_\Theta +\| \widehat{f}^\theta_t- f^\theta_t \|_\Theta ).
 \end{align*}
 According to the proprieties of the function $x \mapsto \log x$, we can show that
 $ \| \log \widehat{f}^\theta_t-\log f^\theta_t\|_\Theta \leq \frac{1}{\underline{c}}\| \widehat{f}^\theta_t- f^\theta_t \|_\Theta$.
Moreover, according to (\textbf{A}$_0 (\Theta)$), we have
\begin{align}
  \| \widehat{f}^\theta_t- f^\theta_t \|_\Theta
   \nonumber   &=  \| f ( Y_{t-1},\cdots,Y_{1},0,\cdots;\theta) - f ( Y_{t-1},\cdots,Y_{1},Y_{0},Y_{-1},\cdots ;\theta) \|_\Theta \\
     \label{proof_lem0_eq1} & \leq \sum\limits_{k\geq t} \alpha^{(0)}_{k} Y_{t-k}.
\end{align}
Hence,
 %
 \[ \frac{1}{\kappa_n} \| \widehat{L}_n(\theta)-{L}_{n}( \theta) \|_\Theta
  \leq   \frac{1}{\kappa_n} \sum_{t=1}^{n}   \Big[ \Big( \frac{Y_t}{\underline{c}}+1 \Big) \| \widehat{f}^\theta_t- f^\theta_t \|_\Theta \Big] 
 \leq   \frac{1}{\kappa_n} \sum\limits_{ \ell=1 }^{n}  \sum\limits_{k \geq \ell } \alpha^{(0)}_{k}  \Big[ \Big( \frac{ Y_{\ell} }{\underline{c}}+1 \Big) Y_{\ell -k} \Big]. \]
 %
 %
 By Corollary 1 of  Kounias and Weng (1969), it suffices to show that
\begin{equation}\label{Kounias_lem0}
 \sum_{\ell \geq 1} \frac{1}{\kappa_\ell} \E \Big[ \sum\limits_{k \geq \ell } \alpha^{(0)}_{k} \Big( \frac{ Y_{\ell} }{\underline{c}}+1\Big) Y_{\ell-k}   \Big]  < \infty.  
 \end{equation}
By H\"{o}lder's inequality and the stationary assumptions, for any $\ell \geq 1$, $k \geq \ell$, it holds that (see (\ref{moment}))
 \[ \E  \Big[ \Big( \frac{ Y_{\ell} }{\underline{c}}+1 \Big) Y_{\ell -k} \Big]  \leq
  \Big( \E\Big[ \Big( \frac{ Y_{\ell} }{\underline{c}}+1 \Big)^2 \Big]  \Big)^{1/2}   \times  (\E Y^2_{\ell-k} )^{1/2}
  = C < \infty .\]
  Hence,
 \[
   \sum_{\ell \geq 1} \frac{1}{\kappa_\ell} \E \Big[ \sum\limits_{k \geq \ell } \alpha^{(0)}_{k} \Big( \frac{ Y_{\ell}     }{\underline{c}}+1\Big) Y_{\ell-k}   \Big]   
   \leq C \sum_{\ell \geq 1}  \frac{1 }{ \kappa_{\ell }}   \sum\limits_{k \geq \ell }   \alpha^{(0)}_{k} < \infty  ,
 \]
where the last equation follows from the condition (\ref{cond_th1}) on the regularization parameter. 
  Hence, (\ref{Kounias_lem0}) holds and thus (\ref{eq_lem0}) follows. 
$~~~~~~~~~~~~~~~~~~~~~~~~~~~~~~~~~~~~~~~~~~~~~~~~~~~~~~~~~~~~~~~~~~~~~~~~~~~~~~~~~~~~~~~~~~~~~~~~~~~~~~~~~~~~~~~~~~~~~~~~ \blacksquare$\\~\\
The following lemma will be useful in the sequel. 
 \begin{lem}\label{lem1}
 Suppose that the assumptions of Theorem \ref{Th1} hold and if a model $m \in \mathcal{M}$ is such that $\theta^* \in \Theta(m)$.   Then
 \[\frac{1}{\kappa_n} \left|\widehat L_n(\widehat  \theta(m)) -\widehat L_n(\widehat{\theta}(m^*)) \right|=o_P(1).\]
\end{lem}

\emph{\bf Proof of Lemma \ref{lem1}}\\
We have
\begin{align*}
\frac{1}{\kappa_n} \left|\widehat L_n(\widehat  \theta(m)) -\widehat L_n(\widehat{\theta}(m^*)) \right|
&= \frac{1}{\kappa_n} \left|\widehat L_n(\widehat  \theta(m)) - L_n(\widehat  \theta(m))+ L_n(\widehat  \theta(m))- L_n(\widehat{\theta}(m^*))+ L_n(\widehat{\theta}(m^*))- \widehat L_n(\widehat{\theta}(m^*)) \right| \\
&\leq \frac{2}{\kappa_n} \left\|\widehat L_n(\theta) - L_n(\theta)  \right\|_\Theta + \frac{1}{\kappa_n} \left| L_n(\widehat  \theta(m))-  L_n(\widehat{\theta}(m^*))\right|\\
&\leq \frac{2}{\kappa_n} \left\|\widehat L_n(\theta) - L_n(\theta)  \right\|_\Theta + \frac{1}{\kappa_n} \left| L_n(\widehat  \theta(m))-  L_n(\theta^*)\right| + \frac{1}{\kappa_n} \left| L_n(\widehat  \theta(m^*))-  L_n(\theta^*)\right| \\
&\leq \frac{2}{\kappa_n} \left\|\widehat L_n(\theta) - L_n(\theta)  \right\|_\Theta + \underset{ \underset{ \theta^* \in \Theta(m) }{ m \in \mathcal{M} } }{\sup}\Big[ \frac{2}{\kappa_n} \left| L_n(\widehat  \theta(m))-  L_n(\theta^*)\right| \Big] 
\end{align*}
Since $\frac{1}{\kappa_n} \left\|\widehat L_n(\theta) - L_n(\theta)  \right\|_\Theta \limitepsn 0$ (from Lemma \ref{lem0}) and $\mathcal{M}$ is a finite collection, it suffices to show that for any $m \in \mathcal{M}$ such that $\theta^* \in \Theta(m)$, 
\begin{equation}\label{eq1_lem1}
\frac{1}{\kappa_n} \left| L_n(\widehat  \theta(m))-  L_n(\theta^*)\right|= o_P(1).
\end{equation}
Let $m \in \mathcal{M}$ with $\theta^* \in \Theta(m)$. By Applying the second order Taylor expansion of $L_n$ around $\widehat \theta(m)$ for $n$ sufficiently large such that $\bar \theta(m) \in \Theta(m)$ which are between $\widehat \theta(m)$ and $\theta^*$, we get
\begin{equation}\label{eq2_lem1}
\frac{1}{\kappa_n} \left( L_n(\widehat  \theta(m))-  L_n(\theta^*)\right)
= 
\frac{1}{\kappa_n}(\widehat \theta(m) -\theta^*) \frac{\partial L_n(\widehat \theta(m))}{\partial \theta}
      + \frac{1}{2\kappa_n} (\widehat \theta(m) -\theta^*)' \frac{\partial^2 L_n(\bar \theta(m))}{\partial \theta \partial \theta'} (\widehat \theta(m) -\theta^*).
\end{equation}
Remark that
\[
\frac{1}{\kappa_n} (\widehat \theta(m) -\theta^*) \frac{\partial L_n(\widehat \theta(m))}{\partial \theta}
    = 
    \frac{1}{\kappa_n}\sqrt{n}(\widehat \theta(m) -\theta^*) \frac{1}{\sqrt{n}}\frac{\partial L_n(\widehat \theta(m))}{\partial \theta}.   
\]
Since $\theta^* \in \Theta(m)$, from \cite{Francq2016}, it holds that $\sqrt{n}(\widehat \theta(m) -\theta^*)=O_p(1)$ and $\frac{1}{\sqrt{n}}\frac{\partial L_n(\widehat \theta(m))}{\partial \theta}=o_P(1)$. 
Hence, 
\begin{equation}\label{eq3_lem1}
\frac{1}{\kappa_n} (\widehat \theta(m) -\theta^*) \frac{\partial L_n(\widehat \theta(m))}{\partial \theta}
    = 
    o_P(1).   
\end{equation} 
Moreover, from \cite{Francq2016}, we get

\[
\sqrt{n}\left(\widehat{\theta}(m)-\theta^*\right) \limiteloin \mathcal{A}_{\theta^*} \equiv  \mathcal{N}(0,J^{-1}(\theta^*) I(\theta^*) J^{-1}(\theta^*)) ~\text{ and } ~
\Big(-\frac{1}{n}\frac{\partial^2 L_n (\bar \theta(m))}{\partial \theta \partial \theta'}\Big) \limitepsn J (\theta^*),
\]
where $J (\theta^*)$ is positive definite.
Hence,
\begin{align*}
 \frac{1}{2} (\widehat \theta(m) -\theta^*)' \frac{\partial^2 L_n(\bar \theta(m))}{\partial \theta \partial \theta'} (\widehat \theta(m) -\theta^*)
 &= \frac{1}{2} \sqrt{n} (\widehat \theta(m) -\theta^*)' \frac{1}{n}\frac{\partial^2 L_n(\bar \theta(m))}{\partial \theta \partial \theta'} \sqrt{n}(\widehat \theta(m) -\theta^*)\\
 &= -\frac{1}{2} \sqrt{n} (\widehat \theta(m) -\theta^*)' (J (\theta^*)+o_p(1)) \sqrt{n}(\widehat \theta(m) -\theta^*)\\
 & \limiteloin -\frac{1}{2} \mathcal{A}'_{\theta^*} J (\theta^*) \mathcal{A}_{\theta^*}.
\end{align*}
%
Therefore, it follows that
\[
(\widehat \theta(m) -\theta^*)' \frac{\partial^2 L_n(\bar \theta(m))}{\partial \theta \partial \theta'} (\widehat \theta(m) -\theta^*)=O_P(1).
\]
Thus, 
\begin{equation}\label{eq4_lem1}
\frac{1}{\kappa_n}(\widehat \theta(m) -\theta^*)' \frac{\partial^2 L_n(\bar \theta(m))}{\partial \theta \partial \theta'} (\widehat \theta(m) -\theta^*)=o_P(1).
\end{equation}
Hence, (\ref{eq1_lem1}) follows from (\ref{eq2_lem1}), (\ref{eq3_lem1}) and (\ref{eq4_lem1}); 
that ends the proof of Lemma \ref{lem1}.

$~~~~~~~~~~~~~~~~~~~~~~~~~~~~~~~~~~~~~~~~~~~~~~~~~~~~~~~~~~~~~~~~~~~~~~~~~~~~~~~~~~~~~~~~~~~~~~~~~~~~~~~~~~~~~~~~~~~~~~~~ \blacksquare$\\

  \medskip
  
 \medskip
 
\noindent Now, let us use the above lemmas to prove Theorem \ref{Th1}.
We proceed as in Bardet {\it et al.} (2019). 
~ \\
 \textbf{(1.)} Firstly, let us prove that $\mathbb{P}(\widehat m_n =m^*) \limiten 1$.
 Remark that
 \[ \mathbb{P}(\widehat m_n =m^*) = 1- \mathbb{P}(m^* \subset \widehat m_n ) - \mathbb{P}(m^* \nsubseteq \widehat m_n ) .\]
 Therefore, it suffices to show that $\mathbb{P}(m^* \subset \widehat m_n ), ~  \mathbb{P}(m^* \nsubseteq \widehat m_n ) \limiten 0$.
 %
 %
 ~ \\
 \textbf{(i)} Note that, if $m^* \subset \widehat m_n$, then there exists a model $m \in \mathcal M$ with $m^* \subset m$ such that $\widehat m_n =m$. Hence,
 \[ \mathbb{P}(m^* \subset \widehat m_n ) \leq \mathbb{P}\big( \underset{   \underset{m\supset m^*}{m \in \mathcal M} }{\cup} \{\widehat m_n =m  \}  \big) \leq \sum_{\underset{m\supset m^*}{m \in \mathcal M}} \mathbb{P}(\widehat m_n =m) .\] 
 Since $\mathcal{M}$ is finite, this item is achieved if we prove that, for any $m \in \mathcal M$ such that $m^* \subset m$, $\mathbb{P}(\widehat m_n =m) \limiten 0$. 
 Let $m \in \mathcal M$ such that $m^* \subset m$. From (\ref{Estim_m}), we have
	\begin{align*}
\mathbb{P}(\widehat m_n =m) &\leq\mathbb{P}\big(\widehat{C}_n(m) \leq \widehat{C}_n(m^*) \big) \\
& \leq\mathbb{P}\Big( -2 \widehat L_n(\widehat \theta(m)) + \kappa_n |m|) \leq -2 \widehat L_n(\widehat \theta(m^*)) + \kappa_n |m^*|  \Big) \\
	&\leq\mathbb{P}\Big(-2\left( \widehat L_n(\widehat  \theta(m))-L_n(\widehat{\theta}(m^*)) \right)  \leq \kappa_n (|m^*|-|m|) \Big)\\
	& \leq \mathbb{P}\Big(\frac{1}{\kappa_n} \big(\widehat L_n(\widehat  \theta(m)) -\widehat L_n(\widehat{\theta}(m^*)) \big)\geq \frac{|m|-|m^*|}{2}\Big)\\
	& ~~\limiten 0, \text{ (according to Lemma}~ \ref{lem1} \text{ and because } |m|> |m^*| \text{)}.
	\end{align*}
 
 \medskip
 \noindent \textbf{(ii)} Similarly as above, let $m \in \mathcal M$ such as $m^* \nsubseteq m$, we are going to show that $\mathbb{P}(\widehat m_n =m) \limiten 0$; which will complete the first part of the proof.  \\
We have 
\begin{align}\label{eq1_th1}
\nonumber \widehat L_n(\widehat  \theta(m^*))-\widehat L_n(\widehat  \theta(m)) &= \left(\widehat L_n(\widehat  \theta(m^*))  -L_n(\widehat  \theta(m^*)) \right) -\left(\widehat L_n(\widehat  \theta(m))  -L_n(\widehat  \theta(m)) \right)\\
&~~~~~~~~~~~~~~~~~~~~~~~~~~~~~~~~~~~~~~~~~~~
+ \left( L_n(\widehat  \theta(m^*))  -L_n(\widehat  \theta(m)) \right).
\end{align}
We can see that by virtue of Lemma \ref{lem0}, the first and the second term of the right part of (\ref{eq1_th1}) are equal to $o(\kappa_n)$.
Moreover, from the proof of Proposition \ref{prop1}, we get
\begin{align}\label{eq2_th1}
\nonumber \frac{1}{n} \big(L_n(\widehat  \theta(m^*))  -L_n(\widehat  \theta(m)) \big) &=  L(\widehat  \theta(m^*)) - L(\widehat  \theta(m))  + o(1) \\
&= L(\widehat  \theta(m^*)) - L( \theta^*) - \big(L(\widehat  \theta(m)) - L(\theta^*(m)) \big) + \big( L( \theta^*) - L(\theta^*(m)) \big) +  o(1) 
\end{align}
where  $L(\theta) = \E[\ell_1(\theta)]$ is defined for all $\theta \in \Theta$.
%
%
According to the consistency of $\widehat  \theta(m^*)$ and $\widehat  \theta(m)$ (see Proposition \ref{prop1}), it holds that
\[ L(\widehat  \theta(m^*)) - L( \theta^*)  = o(1) \text{ and }  L(\widehat  \theta(m)) - L(\theta^*(m)) = o(1) . \] 
Hence, (\ref{eq1_th1}) and (\ref{eq2_th1}) implies 
\begin{equation}\label{proof_th1_approx_theta_start_theta}
  \frac{1}{n} \big( \widehat L_n(\widehat  \theta(m^*))-\widehat L_n(\widehat  \theta(m)) \big) = \frac{1}{n} \big(L_n(\widehat  \theta(m^*))  -L_n(\widehat  \theta(m)) \big) + o(\frac{\kappa_n}{n}) =  L( \theta^*) - L(\theta^*(m)) +  o(1)  .
\end{equation}
Therefore, it comes that
\begin{align}\label{proof_th1_approx_Cm_Cm_start}
 \nonumber \frac{1}{n} \big( \widehat{C}_n( m) - \widehat{C}_n(m^*) \big) &=  
   \frac{2}{n} \big( \widehat L_n(\widehat  \theta(m^*))-\widehat L_n(\widehat  \theta(m)) \big) + \frac{\kappa_n}{n}(|m| - |m^*|) \\ 
    &= 2 \big( L( \theta^*) - L(\theta^*(m)) \big) + \frac{\kappa_n}{n}(|m| - |m^*|) + o(1). 
\end{align}  
 Note that, we can go along the same lines as in proof of Proposition \ref{prop1} to show that the function $L : \Theta \rightarrow \R$, $\theta \mapsto L(\theta) = \E[\ell_1(\theta)]$ has a unique maximum at $\theta^*$ and for all model $m \in \mathcal{M}$, $\theta^*(m) = \theta^*$ when $m \supseteq m^*$.
 Thus, since  $\theta^* \notin \Theta(m)$ (because $m^* \nsubseteq m$), we have $L( \theta^*) - L(\theta^*(m)) >0$.
Hence, according to $\kappa_n = o(n)$ and the fact that $|m|$ and $|m^*|$ are finite, (\ref{proof_th1_approx_Cm_Cm_start}) implies $ \widehat{C}_n( m) > \widehat{C}_n(m^*) $ a.s. for $n$ large enough. Thus,
\[ \mathbb{P}(\widehat m_n =m) \leq\mathbb{P}\big(\widehat{C}_n(m) \leq \widehat{C}_n(m^*) \big)  \limiten 0 .\]

\noindent \textbf{(2.)} The next lines show the second part which is about the consistency of $\widehat \theta(\widehat m_n)$.
Let $\epsilon >0$. We have 
\begin{equation}\label{proof_th1_proba_conv}
 \mathbb{P}( \| \widehat{\theta}(\widehat{m}) - \theta^* \| > \epsilon) =  \mathbb{P}( \| \widehat{\theta}(\widehat{m}) - \theta^* \| > \epsilon | \widehat m_n =m^*) \mathbb{P}(\widehat m_n =m^*) + \mathbb{P}( \| \widehat{\theta}(\widehat{m}) - \theta^* \| > \epsilon | \widehat m_n \neq m^*) \mathbb{P}(\widehat m_n \neq m^*) 
\end{equation}
According to the first part established above, $\mathbb{P}(\widehat m_n =m^*) \limiten 1$ and $ \mathbb{P}(\widehat m_n \neq m^*) \limiten 0$. Thus, the first term of the right hand side of (\ref{proof_th1_proba_conv}) converges to 0 (from the strong consistency of $\widehat{\theta}(\widehat{m}^*)$, see Theorem 2.1 of \cite{Francq2016}) and the second term also converges to 0. This completes the proof of Theorem \ref{Th1}.
 
 $~~~~~~~~~~~~~~~~~~~~~~~~~~~~~~~~~~~~~~~~~~~~~~~~~~~~~~~~~~~~~~~~~~~~~~~~~~~~~~~~~~~~~~~~~~~~~~~~~~~~~~~~~~~~~~~~~~~~~~~~~~~~~~~~~~~~~~~~~~~~~~~~ \blacksquare$

 \subsection{Proof of Theorem \ref{Th2}}
 
 For $x=(x_i)_{1\leq i \leq d} \in \R^d$, set $F_n(x)= \mathbb{P} \big(\underset{1\leq i \leq d}{\bigcap} \sqrt{n}( \widehat \theta(\widehat m_n)-\theta^*)_i \leq x_i \big)$.\\
We have
\begin{align*}
F_n(x) &=
 \mathbb{P} \big(\underset{1\leq i \leq d}{\bigcap} \sqrt{n}( \widehat \theta(\widehat m_n)-m^*)_i \leq x_i \, \big | \,  \widehat m_n=m^* \big)\mathbb{P}(\widehat m_n=m^*)\\
&~~~~~~~~~~~~~~~~~~~~+
 \mathbb{P} \big(\underset{1\leq i \leq d}{\bigcap} \sqrt{n}( \widehat \theta(\widehat m_n)-m^*)_i \leq x_i \, \big | \, \widehat m_n \neq m^* \big)\mathbb{P}(\widehat m_n \neq m^*).
\end{align*}
According to Theorem \ref{Th1},  $\mathbb{P}(\widehat m_n =m^*) \limiten 1$ and $\mathbb{P}(\widehat m_n \neq m^*) \limiten 0$.  Therefore, the second term in the right side of the previous equality converges to zero. The first term can be written as
\begin{align*}
 &\mathbb{P} \big(\underset{1\leq i \leq d}{\bigcap} \sqrt{n}( \widehat \theta(\widehat m_n)-\theta^*)_i \leq x_i \, \big | \,  \widehat m_n=m^* \big)\\
 & ~~~~~~~~~~~~~~~~~~~~~~~~~~=
 \mathbb{P} \Big(\big\{\underset{i \in m^*}{\bigcap} \sqrt{n}( \widehat \theta(m^*)-\theta^*)_i \leq x_i\big\}
 \bigcap
\big\{\underset{i \notin m^*}{\bigcap} \sqrt{n}( \widehat \theta(m^*)-\theta^*)_i \leq x_i\big\}\Big).
\end{align*}
Since $\theta(m^*) \in \Theta(m^*)$, $((\widehat \theta(m^*))_{i})_{i \notin m^*}=(\theta^*_{i})_{i \notin m^*}=0$ and for $(x_i)_{i\notin m^*}$ a family of non-negative real numbers, we have
\begin{align*}
&\mathbb{P} \Big(\big\{\underset{i \in m^*}{\bigcap} \sqrt{n}( \widehat \theta(m^*)-\theta^*)_i \leq x_i\big\}
 \bigcap
\big\{\underset{i \notin m^*}{\bigcap} \sqrt{n}( \widehat \theta(m^*)-\theta^*)_i \leq x_i\big\}\Big)\\
& ~~~~~~~~~~~~~~~~~~~~~~~~~~=
\mathbb{P} \Big(\underset{i \in m^*}{\bigcap} \sqrt{n}( \widehat \theta(m^*)-\theta^*)_i \leq x_i\Big)\\
& ~~~~~~~~~~~~~~~~~~~~~~~~~~\underset{n \rightarrow \infty }{\longrightarrow} 
\mathbb{P} \Big(\big(\Sigma_{\theta^*}\big)^{-1/2}Z\leq (x_i)_{i \in m^*}\Big),
\end{align*}
where $Z$ is the standard Gaussian random vector in $\R^{|m^*|}$ from the central limit given in (\ref{Normal_Theta}); which
completes the proof of Theorem \ref{Th2}.
\\
$~~~~~~~~~~~~~~~~~~~~~~~~~~~~~~~~~~~~~~~~~~~~~~~~~~~~~~~~~~~~~~~~~~~~~~~~~~~~~~~~~~~~~~~~~~~~~~~~~~~~~~~~~~~~~~~~~~~~~~~~~~~~~~~~~~~~~~~~~~~~~~~~ \blacksquare$

\subsection{Proof of Proposition \ref{prop_non_lin}} 
 Recall that
\begin{equation}\label{eq1_prop2}
\E(Y^2_t)=\E\left(\E(Y^2_t|\mathcal{F}_{t-1})\right) = \E\left( Var(Y_t|\mathcal{F}_{t-1}) + (\E(Y_t|\mathcal{F}_{t-1}))^2\right).
\end{equation}
Moreover,
\begin{align*}
Var(Y_t|\mathcal{F}_{t-1})&=r\frac{(1-p_{t})}{p^2_{t}} ~~(\text{because } ~Y_{t}|\mathcal{F}_{t-1} \sim NB(r,p_{t}))\\
            &=\lambda_t+\frac{\lambda^2_t}{r}.
\end{align*}
From (\ref{eq1_prop2}), we deduce that 
\begin{equation}\label{eq2_prop2}
\E(Y^2_t)=\E\big( \lambda_t+\frac{\lambda^2_t}{r} + \lambda^2_t \big)=\E( \lambda_t) +( 1+1/r) E(\lambda^2_t).
\end{equation}
Thus, $\E(Y^2_t)<\infty$ if and only if $\E(\lambda^2_t)<\infty$.

\medskip

\noindent From (\ref{nonlinear_lip}), we can get 
 \begin{equation*}
\lambda_{t} = f_{\theta^*}(Y_{t-1},\lambda_{t-1}) \leq a Y_{t-1}+ b\lambda_{t-1} + f_{\theta^*}(0,0).
 \end{equation*}
 Hence, it follows that
 \begin{align}\label{eq3_prop2}
 \E(\lambda^2_{t}) & \leq \E\bigg( \big(a Y_{t-1}+ b\lambda_{t-1} + f_{\theta^*}(0,0)\big)^2\bigg) \nonumber\\
   & \leq
   \E\big( (a Y_{t-1}+ b\lambda_{t-1})^2\big) +2f_{\theta^*}(0,0) \E(a Y_{t-1}+ b\lambda_{t-1}) + (f_{\theta^*}(0,0))^2 \nonumber\\
    &\leq 
    a^2 \E(Y^2_{t-1})+2ab \E(Y_{t-1}\lambda_{t-1}) + b^2\E(\lambda^2_{t-1})+C \nonumber\\
    &\leq
    a^2 \E \big(\lambda_{t-1} +( 1+1/r)\lambda^2_{t-1} \big)+2ab \E(Y_{t-1}\lambda_{t-1}) + b^2\E(\lambda^2_{t-1})+C
    ~~(\text{from relation }  (\ref{eq2_prop2})) .
 \end{align}
 Remark that 
\begin{equation*}
\E(Y_{t-1}\lambda_{t-1})=\E \big( \E(Y_{t-1}\lambda_{t-1}|\mathcal{F}_{t-2}) \big)=\E(\lambda^2_{t-1}).
\end{equation*}
Thus, from (\ref{eq3_prop2}), we have   
  \begin{align}\label{eq4_prop2}
  \E(\lambda^2_{t}) 
  &\leq 
  a^2 \E \big(\lambda_{t-1} +( 1+1/r)\lambda^2_{t-1} \big)+2ab \E(\lambda^2_{t-1}) + b^2\E(\lambda^2_{t-1})+C \nonumber\\
  &\leq \left( ( 1+1/r)a^2+ 2ab +b^2\right)\E(\lambda^2_{t-1})+C. 
   \end{align}  
Since the process $\{\lambda_{t},~t\geq 1 \}$ is strict stationary, from (\ref{eq4_prop2}), a sufficient condition for that $\lambda_{t}$ has second-order moment is 
 \[
   (1+1/r)a^2+ 2ab +b^2<1;~~i.e~~(a+b)^2+\frac{a^2}{r}<1.
\]  
This achieves the proof of Proposition \ref{prop_non_lin}.
\\
$~~~~~~~~~~~~~~~~~~~~~~~~~~~~~~~~~~~~~~~~~~~~~~~~~~~~~~~~~~~~~~~~~~~~~~~~~~~~~~~~~~~~~~~~~~~~~~~~~~~~~~~~~~~~~~~~~~~~~~~~~~~~~~~~~~~~~~~~~~~~~~~~ \blacksquare$
                     

\end{document}